\documentclass[11pt]{article}

\usepackage{amstex,theorem}

\makeatletter 
\newcommand{\ps@verbit}{%
  \renewcommand{\@oddhead}{%
          \scriptsize
          {Partial resolutions of Hilbert type}
          \hfil\tiny {D. Kaledin, M. Verbitsky, \ \ \ 12 Dec. 1998}}
  \renewcommand{\@evenhead}{\@oddhead}
  \renewcommand{\@oddfoot}{\hfil\thepage\hfil}
  \renewcommand{\@evenfoot}{\@oddfoot}}
\makeatother
 
\long\def\comment#1\endcomment{}

\newtheorem{theorem}{Theorem}[section]
\newtheorem{lemma}[theorem]{Lemma}
\newtheorem{corollary}[theorem]{Corollary}
\newtheorem{prop}[theorem]{Proposition}
\newtheorem{conjecture}[theorem]{Conjecture}

\makeatletter
\begingroup
\gdef\th@upshape{\normalfont
  \def\@begintheorem##1##2{%
        \item[\hskip\labelsep \theorem@headerfont ##1\ ##2]}%
\def\@opargbegintheorem##1##2##3{%
   \item[\hskip\labelsep \theorem@headerfont ##1\ ##2\ (##3)]}}
\endgroup
\makeatother

\theoremstyle{upshape}

\newtheorem{defn}[theorem]{Definition}
\newtheorem{remark}[theorem]{Remark}

\makeatletter
\@addtoreset{equation}{section}
\makeatother

\renewcommand{\labelenumi}{{\normalfont(\roman{enumi})}}

\newcommand{\proof}[1][Proof. ]{\smallskip\noindent{\em #1}}
\def\endproof{\ensuremath{\square}\par\medskip}

\def\eqref#1{\thetag{\ref{#1}}}

\newcommand{\wt}{\widetilde}
\newcommand{\wh}{\widehat}

\newcommand{\idot}{{\:\raisebox{3pt}{\text{\circle*{1.5}}}}}

\newcommand{\C}{{\Bbb C}}
\newcommand{\R}{{\Bbb R}}
\newcommand{\Z}{{\Bbb Z}}
\newcommand{\Q}{{\Bbb Q}}
\newcommand{\h}{{\Bbb H}}

\newcommand{\T}{{\cal T}}
\newcommand{\N}{{\cal N}}
\newcommand{\ZZ}{{\cal Z}}
\newcommand{\calo}{{\cal O}}
\newcommand{\M}{{\cal M}}
\newcommand{\HH}{{\cal H}}

\newcommand{\TT}{{\frak t}}
\newcommand{\n}{{\frak n}}
\newcommand{\U}{{\frak u}}
\newcommand{\su}{{\frak s}{\frak u}}

\renewcommand{\phi}{\varphi}
\newcommand{\arrow}{{\:\longrightarrow\:}}

\newcommand{\Aut}{\operatorname{Aut}}
\newcommand{\End}{\operatorname{End}}
\newcommand{\Norm}{\operatorname{Norm}}
\newcommand{\Card}{\operatorname{Card}}
\newcommand{\Sing}{\operatorname{Sing}}
\newcommand{\Spec}{\operatorname{Spec}}

\newcommand{\Diag}{\operatorname{\sf Diag}}
\renewcommand{\dim}{\operatorname{\sf dim}}

\title{Partial resolutions of Hilbert type, Dynkin diagrams 
and generalized Kummer varieties}

\author{D. Kaledin, M. Verbitsky}

\begin{document}

\maketitle

{\small 
\hspace{0.17\linewidth}
\begin{minipage}[t]{0.7\linewidth}
We study the partial resolutions of singularities related to Hilbert
schemes of points on an affine space.  Consider a quotient of a
vector space $V$ by an action of a finite group $G$ of linear
transforms. Under some additional assumptions, we prove that the
partial desingularization of Hilbert type is smooth only if the
action of $G$ is generated by complex reflections. This is used to
study the subvarieties of a Hilbert scheme of a complex torus. We
show that any subvariety of a generic deformation of a Hilbert
scheme of a torus is birational to a quotient of another torus by a
Weyl group action. In Appendix, we produce counterexamples to a
false theorem stated in our preprint math.AG/9801038.
\end{minipage}
}

\tableofcontents


\section{Introduction}


Let $M$ be a compact hyperk\"ahler manifold. A hyperk\"ahler manifold
is equipped with an action of quaternions in its tangent bundle.
This action induces a set of complex structures (so-called induced
complex structures) on $M$ (see Subsection \ref{_Prelimi_Subsection_}). 
A closed subvariety $X\subset M$
is called {\bf trianalytic} if $X$ is complex analytic with respect
to all induced complex structures
(Definition \ref{_triana_Definition_}). Trianalytic subvarieties
were a subject of a long study. Most importantly, take a generic
induced complex structure $I$ on $M$. Then all closed complex
subvarieties of $(M, I)$ are trianalytic (\cite{v.trian}). 
Moreover,  a trianalytic subvariety can be canonically
desingularized (\cite{_Verbitsky:hypercomple_}),
and this desingularization is hyperk\"ahler.

\hfill

Two series of compact hyperk\"ahler manifolds are studied: the
Hilbert schemes of K3 surfaces and the so-called generalized Kummer
varieties (see Definition \ref{_generali_Kummer_Definition_}).  As
it was shown in \cite{_Verbitsky:hypercomple_}, for a generic
hyperk\"ahler structure on a Hilbert scheme of a K3 surface, there
are no trianalytic subvarieties (except points).  The proof of this
result is technically quite difficult. In \cite{kv}, a similar
result was ``found'' of a generalized Kummer variety. The proof was
elementary, but, unfortunately, false. Indeed, there are
counterexamples (see Appendix), i. e. trianalytic subvarieties of
any deformation of a generalized Kummer variety.

These counterexamples were found as a byproduct of attempts
to prove the conjecture of V. Ginzburg (\cite{_Bezruk_Ginzbu_}),
which is stated roughly as follows.

\begin{conjecture} \label{_Ginz_Conjecture_}
Let $G$ be a group generated by real reflections
acting on a complex vector space $V$ and preserving a rational structure
(more specifically, the Weyl group of some semisimple Lie group
acting on its Cartan algebra). Consider the space
$V\oplus V$ equipped with a diagonal action of $G$, and the
corresponding quotient variety
$X:= (V\oplus V)/G$. Then $X$ can be 
naturally desingularized, and this desingularization
is holomorphically symplectic and admits a hyperk\"ahler 
structure.
\end{conjecture}

A counterexample to \cite{kv} was produced
by A. Kuznetsov in an attempt to prove
Conjecture \ref{_Ginz_Conjecture_} for $G$ a Weyl group
corresponding to the Dynkin diagrams $C_n$.

\hfill

It turns out that trianalytic subvarieties of 
generalized Kummer varieties are deeply related
to Dynkin diagrams. Their relation is the main topic
of this article.

\begin{theorem} \label{_Dynki_triana_Theorem_}
Let $X\subset K^{[n]}$ be a trianalytic subvariety of a generalized
Kummer variety associated with a generic complex torus of 
dimension 2. Then 
\begin{itemize}
\item $X$ is birational to a quotient of 
a torus $T_1$ by an action of a Weyl group $W$ associated with some
reductive Lie group $L$. 
\item The torus $T_1$ is isogeneous to
$T^k$, and $G$ fixes the zero of $T_1$. 
\item The tangent space $T_0 T_1$ at $0 \in T_1$ is idenitied with
${\frak k}\oplus {\frak k}$, where ${\frak k}$ is the Cartan algebra
of $L$. This idenitification is compatible with the action of $W$.
\end{itemize}
\end{theorem}

For a more precise statement of Theorem \ref{_Dynki_triana_Theorem_}
and its proof, see Theorem \ref{trian}. 

\hfill

The idea of a proof of Theorem \ref{_Dynki_triana_Theorem_}
is based on considering a special type of partial resolutions
of singularities, called {\bf Hilbert-type partial resolutions}
(Definition \ref{hilbert.type}). For an introduction to
Hilbert-type partial resolutions, see 
Section \ref{_part_res_simplifi_Section_}.

\hfill

\begin{itemize}

\item In Section \ref{_part_res_simplifi_Section_}, we present a
less general version of the formalism of partial resolutions of
Hilbert type.  In this generality, the statement and the proof of
our results are elementary. None of the results of Section
\ref{_part_res_simplifi_Section_} is used further in this paper.
Section \ref{_part_res_simplifi_Section_} is supposed to be an
elementary introduction to the theory of partial resolutions of
Hilbert type.

\item In Section \ref{_part_res_full_Section_},
we define the partial resolution of Hilbert type
in full generality. We fix the notation and give some preliminary
definitions of geometry of Hilbert schemes. 

\item In Section \ref{linear}, we state and prove results of
geometry of Hilbert-type partial resolutions. We consider
a Hilbert-type partial resolution $X= V^{[\phi]}$ associated with a linear map
$\phi:\; V \arrow W^k$. We show that $X$ is smooth only if 
the normalizer $G:= Norm_{S_k} V$ is a group generated by
complex reflections. If, in addition, the space $W$ is equipped 
with a generic rational lattice and $V$ is rational, then
the group $G$ is generated by {\it real} reflections. 
We prove that $G$ is the Weyl group of a reductive Lie algebra 
acting on $V$ as on two copies of its Cartan algebra.

\item In Section \ref{erratum}, we apply the developed
(purely algebro-geometric) formalism of partial resolutions
of Hilbert type to questions of hyper\-k\"ahler geometry.
Using results of \cite{kv}, we show that
any trianalytic subvariety of a generalized Kummer variety
is a partial resolution of Hilbert type associated with
Hilbert scheme and embeddings of complex tori. This 
allows us to establish the relationship between
trialanytic subvarieties and Dynkin diagrams.

\item In Appendix (Section \ref{_Appe_Exi_Section_}), we 
prove the existence of counterexamples to our earlier 
(false) statements from \cite{kv}. We show that
any deformation of a generalized Kummer variety
contains non-trivial complex analytic 
and trianalytic subvarieties. These subvarieties
are constructed as fixed points of some canonical involution;
therefore, they are smooth.

\end{itemize}


\section{Partial resolutions of Hilbert type -- a simplified version}
\label{_part_res_simplifi_Section_}\label{filler}


None of the results of this section are used further on in this
paper. Here we sacrifice the generality to give a clear view
of the interplay between the groups generated by reflections 
and the symplectic geometry of the Hilbert schemes.

\hfill

Let $W$ be a complex vector space, and $W^{(k)} = W^k/S_k$ the 
symmetric power of $W$, i. e. a quotient of $W^k$ by the
natural action of the symmetric group of $k$ letters, and 
\[ W^{[k]}\stackrel{\pi}{\arrow} W^{(k)}
\]
the Hilbert scheme of $k$ points on $W$. Consider a monomorphism
\[ V \stackrel \phi \hookrightarrow W^k \]
of complex vector spaces. Assume that $\phi(V)$ intersects non-trivially
with 
\begin{equation}\label{_pairwise_non_eq_Equation_}
{W^k}_0:= \{(x_1\neq x_2\neq ... \neq x_k) \in W^k\}
\end{equation}
of pairwise non-equal $k$-tuples of points.
Let ${W^{(k)}}_0$ be the quotient of ${W^k}_0$ by the
action of $S_k$. Consider the natural map
$\sigma:\; {W^k}_0 \arrow {W^{(k)}}_0$,
and let $V^\sigma_0\subset {W^{(k)}}_0$ be the image
of $V_0:= \phi(V) \cap {W^k}_0$ under the map $\sigma$.
Clearly, $V_0^\sigma = V_0/G$, where
$G\subset S_k$ is the normalizer of $V$ in $S_k$.
The map $\pi:\; W^{[k]}\arrow W^{(k)}$
induces an isomorphism
\[ 
   \pi:\; \pi^{-1}\left( {W^{(k)}}_0\right) \arrow  {W^{(k)}}_0.
\]
Let $V_0^{[\phi]}\subset  W^{[k]}$
be the preimage $\pi^{-1}(V_0^\sigma)\subset  W^{[k]}$,
and $V^{[\phi]}$ be a closure of $V_0^{[\phi]}$ in 
the Hilbert scheme $W^{[k]}$.
Clearly, the map $\pi:\; W^{[k]}\arrow W^{(k)}$
induces a birational isomorphism
$\pi:\; V^{[\phi]}\arrow V/G$. 

\begin{defn}
In the above assumptions, the variety $V^{[\phi]}$ is called
{\bf a partial Hilbert-type resolution of $V/G$}.
\end{defn}

This definition is a weaker form of Definition \ref{hilbert.type}.
The difference is, in Definition \ref{hilbert.type} we don't assume
that $V$ necessarily intersects with the set ${W^k}_0$ of pairwise
distinct $k$-tuples.

\begin{theorem}\label{_Hilbert_type_simpli_Theorem_}
Let $V\hookrightarrow W^k$ be an embedding 
of complex vector spaces. Denote the normalizer
$Norm_{S_k}(V)$ by $G$. Assume that 
the image $\phi(V)$ is $GL(W)$-invariant and
intersects non-trivially with
the set ${W^k}_0$ of pairwise distinct $k$-tuples
(\eqref{_pairwise_non_eq_Equation_}).
Consider the partial Hilbert-type resolution
$V^{[\phi]}\subset W^{[k]}$. Assume that the 
variety $V^{[\phi]}$ is smooth.
Then there exists a $G$-invariant decomposition
\[ V= \oplus V_i, i=1, ..., \dim W, \]
such that 
\begin{enumerate}
\renewcommand{\labelenumi}{{\normalfont (\roman{enumi})}}
\item all the spaces $V_i$ are isomorphic as representations of $G$, and
\item the group $G$ acts on $V_i$ by complex reflections
(see Definition \ref{_refle_Definition_} for the definition of
a complex reflection).
\end{enumerate}
\end{theorem}

{\bf Proof:} Consider a commutative subgroup $K\subset GL(W)$
of rank $\dim W$ acting on $W$ by characters. In coordinates,
such action can be written as
\[ (a_1, ... a_n)(x_1, ... x_n) = (a_1 x_1, a_2 x_2, ... a_n x_n).
\]
Let $K_i \subset K$, $i=1, ..., \dim W$ be a codimension 1 subgroup
\[ K_i = (a_1,..., a_{i-1}, 1, a_{i+1}, ..., a_n).
\]
The space of invariants of the action of $K_i$ on $W$ is 1-dimensional,
and we denote it by $W_i$. Let $V_i$ be the space of invariants of
the action of $K_i$ on $\phi(V)$. 
Clearly, $V= \oplus V_i$. Moreover, the action of $GL(W)$ on $\phi(V)$ 
commutes with the action of $G$ on $\phi(V)$, and therefore,
the decomposition $V= \oplus V_i$ is $G$-invariant. Using
an element $\eta\in GL(W)$ which interchanges $K_i$ and $K_j$,
we obtain a $G$-equivariant isomorphism $\eta:\; V_i \arrow V_j$.
This proves Theorem \ref{_Hilbert_type_simpli_Theorem_} (i).
It remains to show that the action of $G$ on $V_i$ (say, on $V_1$) 
is generated by complex reflections. Clearly, $\phi$ maps
$V_1$ to $W_1 ^k$. Denote this map by $\phi_1:\; V_1 \arrow W_1^k$.
Let $V_1^{[\phi_1]}$ be the corresponding Hilbert-type
partial resolution. Consider the natural embedding
$W_1^{[k]}\stackrel t\hookrightarrow W^{[k]}$ obtained from the functoriality
of Hilbert schemes. 

\begin{lemma} 
The image
$t(V_1^{[\phi_1]})$ lies
in $V^{[\phi]}$, and, moreover, $t(V_1^{[\phi_1]})$ coincides
with the set
\[ \{x\in V^{[\phi]}\;\; | \;\; K_1(x) =x\}
\]
of $K_1$-invariant points in $V^{[\phi]}$
\end{lemma}

{\bf Proof:}
Clear from definitions. \endproof

Since $V_1^{[\phi_1]}$ is a fixed point set of a reductive 
algebraic group acting on a complex manifold, this variety is smooth.
On the other hand, the projection $\pi:\; W_1^{[k]}\arrow W_1^{(k)}$
is an isomorphism, because $\dim W_1 =1$. This implies that the
natural projection $\pi :\; V_1^{[\phi_1]}\arrow V_1/G$ is an isomorphism.
Since $V_1^{[\phi_1]}$ is smooth, $V_1/G$ is also smooth.
Now, by \cite{_Bourbaki:Lie4-6_}, Ch. V, \S 5 Theorem 4
(see also Theorem \ref{smooth.quotient}), 
the quotient $V_1/G$ is smooth if and only if
$G$ acts on $V_1$ by complex reflections. This finishes
the proof of Theorem \ref{_Hilbert_type_simpli_Theorem_}.
\endproof


\section{Partial resolutions of Hilbert type}
\label{_part_res_full_Section_}


\subsection{General definitions}

Let $X$ be an irreducible complex variety. To fix the 
terminology, we
introduce the following.

\begin{defn}\label{partial}
An irreducible complex variety $\wt{X}$ equipped with a proper
morphism $f:\wt{X} \mapsto X$ is called a {\em partial resolution}
of the variety $X$ if it is an isomorphism outside of the subset
$\Sing X \subset X$ of singular points in $X$.
\end{defn}

Assume given a smooth complex variety $X$ of dimension $\dim X >
1$. Let $k > 1$ be an integer, and let $X^k = X \times \cdots \times
X$ be the $k$-fold self-product of the variety $X$. Denote by
$X^{(k)} = X^k / S_k$ the quotient of the variety $X^k$ with respect
to the natural action of the symmetric group $S_k$ on $k$
letters. Let $X^{[k]}$ be the Hilbert scheme of $0$-dimensional
subschemes in $X$ of length $k$. Then we have a canonical proper
projection $\pi:X^{[k]} \to X^{(k)}$ which is a partial resolution
in the sense of Definition~\ref{partial}.

In this paper we will study the following situation. Consider a
smooth submanifold $Y \subset X^k$. Denote by $\phi:Y \to X^k$ the
embedding map. Let $G = \Norm (Y) \subset S_k$ be the normalizer
subgroup of $Y$ in the symmetric group $S_k$, that is, the subgroup
of elements $s \in S_k$ such that $s(Y) = Y \subset X^k$. Consider
the quotient variety $Y/G$. The embedding map $\phi:Y
\hookrightarrow X^k$ defines a natural closed embedding $Y/G \to
X^{(k)} = X^k / S_k$.

Let $(Y/G)^{[\phi]} \subset X^{[k]}$ be a closed subvariety such
that the canonical projection $\pi:X^{[k]} \to X^{(k)}$ maps
$(Y/G)^{[\phi]} \subset X^{[k]}$ onto $Y/G \subset X^{(k)}$. 

\begin{defn}\label{hilbert.type}
If the induced map $\pi:(Y/G)^{[\phi]} \to Y/G$ is a partial
resolution in the sense of Definition~\ref{partial}, then it is
called a {\em partial resolution of Hilbert type}.
\end{defn}

Note that a partial resolution of Hilbert type $(Y/G)^{[\phi]}$
depends not only on the complex variety $Y$ equipped with an action
of the group $G$, but also on the embedding $\phi:Y \hookrightarrow
X^k$. Moreover, it depends on the particular choice of the closed
subset $(Y/G)^{[\phi]} \subset X^{[k]}$ lying over $Y/G \subset
X^{(k)}$. If the subset $Y/G \subset X^{(k)}$ intersects the open
subset
$$
X^{(k)}_0 = \left\{(x_1,\ldots,x_k) \in X^{(k)}|x_1 \neq \cdots \neq
x_k\right\} \subset X^{(k)}, 
$$
then this latter choice is unique, since the projection $\pi:X^{[k]}
\to X^{(k)}$ is bijective over $X^{(k)}_0 \subset X^{(k)}$.

Definition~\ref{hilbert.type} in full generality is probably
useless. In this paper we will consider only two particular cases,
where this definition leads to interesting results. These cases are
\begin{enumerate}
\item $X = W$ is a complex vector space, $Y = V \subset W^k$ is a
linear subspace, $\phi:V \hookrightarrow W^k$ is the given
embedding, or, more generally, an affine embedding parallel to the
given one.
\item $X = T$ is a complex torus, $Y \subset T^k$ is a subtorus,
$\phi:Y \hookrightarrow T^k$ is the natural embedding.
\end{enumerate}

Note that Definition~\ref{hilbert.type} for vector spaces
generalizes the notion of a {\em Hilbert resolution} introduced by
Y. Ito and H. Nakajima in \cite{Nakajima.McKay}. In this case, we have $X
= Y = V$ for some vector space $V$. The finite group $G$ acting on
$V$ is given {\em a priori}. The integer $k = \Card(G)$ equals the
number of elements in $G$, and the embedding $V \hookrightarrow V^k
= V^{{\Card(G)}}$ coincides with the coaction $V \to V \otimes \C(G)
= V^{{\Card(G)}}$ on the space $V$ of the coalgebra $\C(G)$ of
functions on the group $G$. It is easy to check that the given group
$G$ coincides with the normalizer subgroup of the subspace $V
\subset V^{{\Card{G}}}$ in the symmetric group on $\Card{G}$
letters.

In this paper we will be interested in another case, namely, the
case where $X$ is a vector space $W$ of dimension $2$, while the
subvariety $Y \subset X^k$ is a vector subspace $V \subset W^k$ of
dimension higher than $2$. Starting with Section~\ref{linear}, $W =
\C^2$ will always denote the standard complex vector space of
dimension $2$.

\subsection{Local charts}

We will also need a local version of the case (i) of
Definition~\ref{hilbert.type}, where a vector space $W=X$ is
replaced with an open neighborhood $U \subset W$ of $0 \in W$, and
$V \subset W^k$ is replaced with $V \cap U^k \subset W^k$. This will
be important, since every point $x \in X^{[k]}$ the Hilbert scheme
$X^{[k]}$ for a general smooth variety $X$ of dimension $n$ has an
open neighborhood which is a product of Hilbert schemes for a
coordinate neighborhood $U \subset \C^n$. We will now describe this
in some detail.

Let $x \in X^{[k]}$ be such a point, and assume that the
corresponding $0$-dimensional subscheme $\ZZ \subset X$ is supported
on a finite subset $\{x_1,\ldots,x_l\} \subset X$ of $l$ distinct
points in $X$. Further, assume that the part of $\ZZ \subset x$
which is supported on $x_i$ has length $a_i$. The numbers
$a_1,\ldots,a_l$ form a partition
$$
\Delta = \{a_1,\ldots,a_l\} 
$$
of the integer $k$ (which we define as a set of numbers
$a_1,\ldots,a_l$ such that $\sum a_l = k$). For every one of the
points $x_i$, choose a coordinate neighborhood $U_i \subset X$ and
an identification $h_i:U \cong U_i$ between $U_i$ and a fixed open
neighborhood $U \subset \C^n$ of $0 \in \C^n$.

Assume that if $j \neq j$ then $U_i \cap U_j = \emptyset$. Consider
the open subset $U_x \subset X^{[k]}$ of subschemes $\ZZ \subset X$
supported on $\bigcup U_i$ in such a way that $\ZZ \cap U_i$ is of
length $a_i$. Then we have
$$
U_x \cong U_1^{[a_1]} \times \cdots \times U_l^{[a_l]},
$$
and the local isomorphisms $h_i:U \cong U_i$ identify $U_x$ with
the product
$$
U^{[\Delta]} = U^{[a_1]} \times \cdots \times U^{[a_l]}.
$$
We will say that the open subset $U^{[\Delta]} \cong U_x \subset
X^{[k]}$ is a {\em local chart} in $X^{[k]}$ near the point $x \in X^{[k]}$. 

Note that under the canonical projection $\pi:X^{[k]} \to X^{(k)}$,
the local chart $U_x \cong U^{[\Delta]}$ is mapped onto an open
subset $U_{\pi(x)} \subset X^{(k)}$ in the symmetric power
$X^{(k)}$. This subset depends only on $\pi(x) \in X^{(k)}$ and not
on the particular choice of $x \in \pi^{-1}(\pi(x))$. Moreover, we
have $U_x = \pi^{-1}(U_{\pi(x)}) \subset X^{[k]}$.  We will say that
the open subset $U_{\pi(x)} \subset X^{(k)}$ is a {\em local chart}
in $X^{(k)}$ near the point $\pi(x) \in X^{(k)}$.

Let $S_\Delta = S_{a_1} \times \cdots \times S_{a_l} \subset S_k$ be
the subgroup in the symmetric group associated to the partition
$\Delta = \{a_1,\ldots,a_l\}$ of the integer $k$. Then the local
chart $U_{\pi(x)} \subset X^{(k)}$ is canonically isomorphic to the
quotient variety $U^{(\Delta)} = U^k/S_{\Delta} = U^{(a_1)} \times
\cdots \times U^{(a_l)}$.

We will use these local charts to study locally partial resolutions
of Hilbert type. To do this, we will need a slight generalization of
Definition~\ref{hilbert.type}. We will formulate it for a general
smooth manifold $X$, but we will use it only for $X = \C^n$ and for
an open neighborhood $X = U \subset \C^n$ of $0 \subset \C^n$.

Fix a partition $\Delta = \{a_1,\ldots,a_l\}$ of the integer $k$.
Let $X^{[\Delta]} = X^{[a_1]} \times \cdots \times X^{[a_l]}$ be the
product of Hilbert schemes of subschemes in $X$ of lengths
$a_1,\ldots,a_l$. The product of the natural maps $X^{[a_i]} \to
X^{(a_i)}$ gives a natural partial resolution $\pi:X^{[\Delta]} \to
X^{(\Delta)}$.

Let $Y \subset X^k$ be a closed submanifold, let $\phi:X
\hookrightarrow X^k$ be the embedding map, and let $G = \Norm Y
\subset S_{\Delta}$ be the normalizer subgroup of $Y$ in the group
$S_{\Delta}$. Then we have a canonical embedding $\phi:Y/G \to
X^{(\Delta)}$.

\begin{defn}\label{ext.hilbert.type}
A closed subvariety $(Y/G)^{[\phi]} \subset X^{[\Delta]}$ is called
an {\em extended partial resolution of Hilbert type} of the quotient
$Y/G \subset X^{(\Delta)}$ iff the canonical projection
$\pi:X^{[\Delta]} \to X^{(\Delta)}$ maps $(Y/G)^{[\phi]} \subset
X^{[\Delta]}$ onto $Y/G \subset X^{(\Delta)}$, and the induced map
$\pi:(Y/G)^{[\phi]} \to Y/G$ is a partial resolution in the sense of
Definition~\ref{partial}.
\end{defn}

This definition gives a local counterpart to
Definition~\ref{hilbert.type} in the following sense.  Let $Y
\subset X^k$ be a smooth subvariety with normalizer subgroup $G =
\Norm Y \subset S_k$, and let $(Y/G)^{[\phi]} \subset X^{[k]}$ be a
partial resolution of Hilbert type in the sense of
Definition~\ref{hilbert.type}. Choose a point $x \in Y$, and let
$G_0 \subset G$ be the subgroup of elements in $G$ fixing the point
$x \in Y$. The point $y \in Y \subset X^k$ defines a point $\wt{x}
\in Y/G \subset X^{(k)}$ in the quotient variety $X^{(k)} =
X^k/S_k$.

Let $U^{(\Delta)} \subset X^{(k)}$ be a local chart in $X^{(k)}$
near the point $\wt{x} \in X^{[k]}$, and let $U^{[\Delta]} =
\pi^{-1}(U^{(\Delta)} \subset X^{[k]}$ be the corresponding local
chart in $X^{[k]}$. Then we have
$$
Y/G \cap U^{(\Delta)} = (Y \cap U^k)/G_0 \subset X^{(k)},
$$
and the intersection $U^{[\Delta]} \cap (Y/G)^{[\phi]} \subset
U^{[\Delta]}$ is an extended partial resolution of Hilbert type of
the quotient $(Y \cap U^k)/G_0$ in the sense of
Definition~\ref{ext.hilbert.type}.

\subsection{Lattices and partial resolutions associated with complex tori} 

\begin{defn}\label{lattice}
A {\em lattice} $L$ in a complex vector space $V$ is a
$\Z$-submodule $L \subset V$ such that the induced map $L \otimes \R
\to V$ is an isomorphism. A {\em rational lattice} $L_\Q \subset V$
is a $\Q$-vector subspace such that $L_\Q \otimes_\Q \R \cong V$.
\end{defn}

For every lattice $L \subset V$ in a complex vector space $V$, the
tensor product $L \otimes \Q \subset V$ is a rational lattice. Two
lattices $L,L' \subset V$ are called isogenic if $L \otimes \Q = L'
\otimes \Q \subset V$. We will say that a subspace $V' \subset V$ in
a complex vector space $V$ {\em compatible} with a lattice $L
\subset V$ in $V$ if the intersection $L' = L \cap V' \subset V'$ is
a lattice in $V'$.

\begin{defn}\label{generic}
A lattice, resp. a rational lattice $L \subset V$ in a complex
vector space $V$ is called {\em generic} if every endomorphism of
the space $V$ preserving $L \subset V$ is a multiplication by an
integer, resp. a rational number.
\end{defn}

Pairs $\langle V, L_\Q \subset V \rangle$ of a complex vector space
$V$ and a rational lattice $L_\Q \subset V$ form a semisimple
$\Q$-linear abelian category. (In fact, it is equivalent to the
category of pure $\Q$-Hodge structures of weight $1$ with only
non-trivial Hodge numbers $h^{1,0}$ and $h^{0,1}$.) Note that if a
rational lattice $L \subset V$ in a vector space $V$ is generic,
then the pair $\langle V, L \otimes \Q \rangle$ is an irreducible
object.

Let $T$ be a complex torus. The torus $T$ is isomorphic to the
quotient $T = W/L$, where $W = \Gamma(T,\T(T))$ is the space of
global holomorphic sections of the holomorphic tangent bundle
$\T(T)$ to the torus $T$, and $L \subset W$ is a lattice in the
complex vector space $W$.

Let $T' \subset T^k$ be a subtorus in the $k$-fold self-product $T^k
= W^k/L^k$. Let $T' = V/L'$, where $V$ is a complex vector space and
$L' \subset V'$ is a lattice in $V'$. The embedding $\phi:T' \to
T^k$ defines a canonical linear embedding $\phi:V' \hookrightarrow
W^k$. The subspace $V' \subset W^k$ is compatible with the lattice
$L^k \subset W^k$, and the lattice $L'$ is the intersection $L' = V'
\cap L^k \subset V'$.

The tangent bundle $\T(T)$ to the torus $T = W/L$ is trivial, and
one can choose a frame consisting of commuting vector fields. For
every point $t \in T$, the associated exponential map defines a
holomorphic covering $W \to T$. Consequently, every sufficitenly
small open neighborhood $U \subset W$ of $0 \subset W$ defines a
canonical coordinate neighborhood $U \cong U_t \subset T$. We will
call this neighborhood {\em flat}. The $k$-fold product $U^k \cong
U_{t_1} \times \cdots \times U_{t_k} \subset T^k$ of flat coordinate
neighborhoods is flat, and the intersection $U^k \cap T'$ with the
subtorus $T' = V/L' \subset T^k$ coincides with the intersection
$U^k \cap V \subset W^k$ with the linear subspace $V \subset
W^k$. The local charts near a point in $T^{(k)}$ or $T^{[k]}$
associated to flat coordinate neighborhoods will also be called
flat.

Let $G = \Norm T' \subset S_k$ be the normalizer subgroup of the
subtorus $T' \subset T^k$ in the symmetric group $S_k$ acting on
$T^k$. Consider the corresponding subvariety $T'/G \subset T^{(k)}$
in the symmetric power $T^{(k)}$. Let $T^{[k]}$ be the Hilbert
scheme of $0$-dimensional subschemes of length $k$ in $T$, and let
$(T'/G)^{[\phi]} \subset T^{[k]}$ be a partial resolution of Hilbert
type of $T'/G \subset T^{(k)}$ in the sense of
Definition~\ref{hilbert.type}.

Consider a point $t \in T' \subset T^k$ in the subtorus $T' \subset
T$, and denote by $G_0 \subset G$ the subgroup of elements in $G =
\Norm T' \subset S_k$ fixing $t \in T'$. Let $\wt{t} \in T'/G
\subset T^{(k)}$ be the associated point in the quotient $t^{(k)}$,
and assume that $U^{(\Delta)} \subset T^{(k)}$ is a flat local
chart in $T^{(k)}$ near the point $\wt{t} \in T^{(k)}$. Then the
intersection $T'/G \cap U^{(\Delta)}$ coincides with the quotient
$U_0/G_0$ for an open neighborhood $U_0 = U^k \cap V \subset V$ of
$0 \in V$, and the intersection $U^{[\Delta]} \cap (T'/G)^{[\phi]}$
provides an extended partial resolution of Hilbert type for this
quotient.


\section{Partial resolutions of Hilbert type and groups generated by
reflections}
\label{linear} 


\begin{defn}\label{_refle_Definition_}
Let $V$ be a complex vector space. Recall that an automorphism $g:V
\to V$ of the vector space $V$ is called a {\em complex reflection}
if it is of finite order and the subspace of invariants $V^g \subset
V$ is of codimension exactly $1$. A complex reflection $g:V \to V$ is
called {\em real} if it is an automorphism of order $2$
(equivalently, if $g$ preserves a real structure on the complex
vector space $V$).
\end{defn}

Finite groups $G \subset \Aut V$ of automorphisms of a complex
vector space $V$ which are generated by complex reflections have
been an object of much study, and there exists a classification of
pairs $\langle V, G \rangle$ of this type. In particular, a subgroup
$G \subset \Aut V$ is generated by {\em real} reflections if and
only if it is a product of Weyl groups associated to Dynkin diagrams
of finite type, assuming that $G$ preserves 
some rational lattice in $V$ 
(see \cite{_Bourbaki:Lie4-6_}, Ch. VI, \S 2, p. 5, Proposition 9).

Let now $W$ be the standard $2$-dimensional complex vector space
equipped with the canonical action of the group $U(2)$. Let $V
\subset W^k$ be a vector subspace, and let $G = \Norm(V) \subset
S_k$ be the normalizer subgroup of the subspace $V \subset W^k$. Let
$\wt{G} \subset \Aut V$ be the quotient of the group $G$ by the
subgroup $G_0 \subset G$ of elements which act trivially on the
vector space $V$.

The standard vector space $W$ carries a natural action $\tau:U(2)
\to \Aut W$ of the group $W(2)$ of unitary $2 \times
2$-matrices. This actions defines $W(2)$-actions on the vector space
$W^k$, on the symmetric power $W^{(k)}$ and on the Hilbert scheme
$W^{[k]}$.

Assume that the subspace $V \subset W^k$ is preserved by the
$U(2)$-action. Choose once and for all a subgroup $U(1) \subset
U(2)$ such that the space of invariants $W^{U(1)} \subset W$ is of
dimension $1$. Then the subspace of invariants $V^{U(1)} \subset V$
is of dimension one half of the dimension of the vector space $V$,
and we have a canonical $U(2)$-equivariant isomorphism $V \cong W
\otimes_\C V^{U(1)}$. (The $U(2)$-action on the right-hand sie is
induced by the standard action $\tau$ on the space $W$.)

Since the action of symmetric group $S_k$ on the space $W^k$
commutes with the $U(2)$-action, the $G$-action on $V \subset W^k$
is also $U(2)$-equivariant. Therefore the group $G$ preserves the
subspace $V^{U(1)} \subset V$. If we equip the space $W$ with the
trivial $G$-action, then the canonical isomorphism $V \cong W
\otimes V^{U(1)}$ is $G$-equivariant.

In this section we prove the following two results. 

\begin{theorem}\label{complex.refl}
Let $(V/G)^{[\phi]} \subset W^{[k]}$ be a partial resolution of
Hilbert type of the quotient $V/G$ in the sense of
Definition~\ref{hilbert.type}. Assume that the closed subvariety
$(V/G)^{[\phi]} \subset W^{[k]}$ is smooth and invariant under the
canonical $U(2)$-action on the Hilbert scheme $W^{[k]}$. Then the
group $\wt{G} \subset \Aut V$ acting on the vector space $V^{U(1)}$
is generated by complex reflections.
\end{theorem}

\begin{remark}
The condition of $U(2)$-invariance imposed on a partial resolution
of Hilbert type $(V/G)^{[\phi]} \subset W^{[k]}$ seems to be quite
non-trivial. Note, however, that it holds automatically when the
vector subspace $V \subset W^k$ is $U(2)$-invariant and intersects
non-trivially with the subset $W^k_0 \subset W^k$ consisting of
$k$-tuples of distinct elements in $W$. This is the simplified case
that we have considered in Section~\ref{filler}.
\end{remark}

\begin{theorem}\label{real.refl}
In the assumptions of Theorem~\ref{complex.refl}, assume
additionally that the subspace $V \subset W^k$ is compatible with a
lattice $L \subset W$ which is generic in the sense of
Definition~\ref{generic}. Then the
subgroup $\wt{G} \subset \Aut V^{U(1)}$ is generated by real
reflections 
\end{theorem}

\begin{remark}
{}From Theorem~\ref{complex.refl} it follows that the action of
$\wt{G}$ is generated by complex reflections. 
\end{remark}

Before we prove Theorem~\ref{complex.refl}, we need to recall a fact
on the $U(2)$-action on the Hilbert scheme $W^{[k]}$. Let $W^{(k)} =
W^k/S_k$ be the symmetric power of the space $W$, and let
$\pi:W^{[k]} \to W^{(k)}$ be the canonical projection. The
projection $\pi$ is compatible with the natural $U(2)$-actions on
$W^{[k]}$ and $W^{(k)}$. Recall that we have fixed a subgroup $U(1)
\subset U(2)$ such that $W^{U(1)} \subset W$ is of dimension
$1$. For every point $u \in W^{(k)}$ let $F_u = \pi^{-1}(u) \subset
W^{[k]}$ be the fiber of the map $\pi$ over the point $u$. (Any such
fiber is a product of so-called {\em punctual Hilbert schemes}.) We
will need the following.

\def\whodhavethought{see also \cite[Section~5]{Nakajima.Hilbert}}

\begin{lemma}[\whodhavethought]\label{nakajima} 
Let $u \in W^{(k)}$ be a point invariant under the action of the
fixed subgroup $U(1) \subset U(2)$. Then the closed subvariety
$F_u^{U(1)} \subset F_u$ of points in the fiber $F_u \subset
W^{[k]}$ invariant under the $U(1)$-action consist of a finite
number of points.
\end{lemma}

\proof Choose a basis $x,y \in W^*$ in the space of linear functions
on the vector space $W$ compatible with the $U(1)$-action, so that
$x$ is $U(1)$-invariant and $y$ satisfies $\lambda \cdot y = \lambda
y$ for every $\lambda \in U(1) \subset \C$. Let $\ZZ \subset W$ be a
$0$-dimensional subscheme of length $k$ in $W$ which is supported at
$0 \subset W$ and invariant under the $U(1)$-action. Then the
$U(1)$-action induces a grading
$$
\calo = \bigoplus_i \calo_i
$$
on the space $\calo = \calo_\ZZ$ of functions on the scheme $\ZZ$, 
defined by
$$
\lambda \cdot f = \lambda^if, \qquad f \in \calo_i, \lambda \in
U(1) \subset \C.
$$
Multiplication by $x$ preserves each of the subspaces $\calo_i =
\calo$. Moreover, since the subscheme $\ZZ \subset W$ is supported
at $0$, multiplication by $x$ is nilpotent. Let $n_i$ be the
smallest positive integer such that $x^{n_i}\calo_i=0$. We obviously
have $n_i \leq \dim_\C \calo_i$.

For every $i$, we have $x^{n_i}y^i\calo = 0$. Indeed, it suffices to
prove that $x^{n_i}y^i \cdot 1 = 0$, where $1 \subset \calo$ is the
unity. But this follows from definition, since $1 \in \calo_0$ and
$y^i \cdot 1 \in \calo_i$. Therefore the ideal
$$
I = \langle x^{n_i}y^i \rangle \subset \C[x,y]
$$
in the polynomial algebra $\C[x,y]$ generated by monomials
$x^{n_i}y^i$ annules the $\C[x,y]$-module $\calo$.

On the other hand, since $n_i \leq \dim_\C\calo_i$, we have
$$
\dim_\C(\C[x,y]/I) = \sum_i n_i \leq \sum_i \dim_\C \calo_i =
\dim_\C \calo = k.
$$
Therefore the surjective map $\C[x,y]/I \to \calo$ defined by $f
\mapsto f \cdot 1$ is in fact bijective, and the subscheme $\ZZ =
\Spec\C[x,y]/I \subset W$ is uniquely defined by the set of natural
numbers $n_i$ satisfying 
$$
\sum_i n_i = k.
$$
There exists only a finite number of such sets. \endproof

\proof[Proof of Theorem~\ref{complex.refl}.]  Let $V_0 = V^{U(1)}
\subset V$ be the subspace of $U(1)$-\-in\-va\-ri\-ant vectors. The
subset $(V/G)^{U(1)} \subset V/G \subset W^{(k)}$ of points
invariant under the $U(1)$-action obviously coincides with $V_0/G
\subset V/G$.

Consider the closed subvariety 
$$
D = \left((V/G)^{[\phi]}\right)^{U(1)} \subset (V/G)^{[\phi]}
$$ 
of $U(1)$-invariant points in the partial resolution $(V/G)^{[\phi]}
\subset W^{[k]}$. Since the variety $(V/G)^{[\phi]}$ is smooth,
while the group $U(1)$ is compact, the subvariety $D \subset
(V/G)^{[\phi]}$ is a union of a finite number of smooth connected
components $D_0,\ldots,D_i$.

The projection $\pi:(V/G)^{[\phi]} \to V/G$ is $U(2)$-equivariant,
therefore it maps the subset $D$ into the subset $V_0/G \subset V/G$
of $U(1)$-invariant points in the quotient $V/G$. Moreover, by
Definition~\ref{partial} the projection $\pi$ is one-to-one over
non-singular points in $V/G$.  Since the generic point of the
subvariety $V_0/G \subset V/G$ is non-singular in $V/G$, there
exists one and only one connected component of the variety $D$, say,
$D_0$, which maps onto $V_0/G \subset V/G$ in such a way that the
map $\pi:D_0 \to V_0/G$ is generically one-to-one.

But by Lemma~\ref{nakajima} the projection $\pi:D_0 \to V_0/G$ is
finite. Since the variety $D_0$ is normal, the finite dominant
projection $\pi:D_0 \to V_0/G$ is the normalization of the variety
$V_0/G$. However, the variety $V_0/G$, being a quotient of a vector
space by a finite group action, is itself normal. Therefore $\pi$ is
an isomorphism between $D_0$ and $V_0/G$. In particular, the
quotient $V_0/G$ is smooth. To finish the proof of
Theorem~\ref{complex.refl}, it remains to invoke the following
classic result.

\begin{theorem}[\cite{_Bourbaki:Lie4-6_}, \ Ch. V, \ \S 5, \ Theorem 4]
\label{smooth.quotient}
The quotient $V/G$ of a {complex} vector space $V$ by a finite
subgroup $G \subset \Aut V$ is smooth if and only if the subgroup $G
\subset \Aut V$ is generated by complex reflections. \endproof
\end{theorem}

\proof[Proof of Theorem~\ref{real.refl}.]
By assumption the intersection $L' = V \cap L^K \subset W^k$ is a
lattice in the vector space $V \subset W^k$. By
Theorem~\ref{complex.refl}, it suffices to prove the following:
\begin{itemize}
\item Any element $g \in G$ which acts as a complex reflection on the
space $V_0 = V^{U(1)}\subset V$ and preserves the lattice $L'
\subset V$ is in fact a real reflection.
\end{itemize}
Given such an element $g$, let $V^g \subset V$ be the subspace of
$g$-invariant vectors in $V$, and let $V_0^g \subset V_0$ be the
subspace of $g$-invariant vectors on $V_0$. Since $V \cong V_0
\otimes W$, we have $V/V^g \cong W \otimes (V_0/V_0^g)$. By
assumption $V_0/V_0^g$ is $1$-dimensional, therefore the quotient
$V/V^g$ is of complex dimension $2$. The element $g$ acts naturally
on the quotient $W \cong V/V^g$, and it suffices to prove that it
acts as an automorphism of order $2$.

Consider the rational lattice $L' \otimes \Q \subset V'$. The
category of complex vector spaces equipped with a rational lattice
is abelian and semisimple.  Since the element $g$ preserves the
lattice $L'$, the quotient $V/V^g$ is equipped with a canonical
quotient rational lattice $L_0 \subset V/V^g$. Moreover, by
assumption the lattice $L \subset W$ is generic in the sense of
Definition~\ref{generic}. Consequently, the object $\langle W, L
\otimes \Q \rangle$ is irreducible. Since $\langle W^k, L^k \otimes
\Q \rangle$ is a sum of several copies of the irreducible object
$\langle W, L \otimes \Q \rangle$, the subobject $\langle V', L'
\otimes \Q \rangle$ is also a sum of several copies of $\langle W, L
\otimes \Q \rangle$, and so is the quotient object $\langle V/V^g,
L_0 \rangle$. Since the vector space $V/V^g$ is $2$-dimensional,
this implies that
$$
\langle V/V^g, L_0 \rangle \cong \langle W, L \otimes \Q \rangle.
$$
Thus the element $g$ is an automorphism of the irreducible object
$\langle W, L \otimes \Q \rangle$. By Definition~\ref{generic} it
must act as multiplication by an invertible rational number, in
other words, by $\pm 1$. Since $g$ is non-trivial, it acts as the
multiplication by $-1$. This finishes the proof. 
\endproof

To simplify notation, we have formulated Theorems~\ref{complex.refl}
and \ref{real.refl} for submanifolds of the Hilbert scheme
$W^{[k]}$.  However, the same results hold for extended partial
resultions of Hilbert type (Definition~\ref{ext.hilbert.type}),
which by definition lie in the product of Hilbert schemes
$W^{[\Delta]}$ associated to a partition $\Delta =
\{a_1,\ldots,a_l\}$ of the integer $k$. To see this, it suffices to
note that the fibers of the canonical projection $\pi:W^{[\Delta]}
\to W^{(\Delta)}$ are products of the fibers of the projections
$\pi:W^{[a_i]} \to W^{(a_i)}, i = 1,\ldots,l$. Therefore
Lemma~\ref{nakajima} is also valid for $W^{[\Delta]}$. Consequently,
the formulations and the proofs of Theorems~\ref{complex.refl} and
\ref{real.refl} carry over to the case of extended resolutions
word-by-word. Moreover, both theorems hold (with the same proofs)
even if we are only given an extended partial resolution for the
quotient $(V \cap U^k)/G$, where $U \subset W$ is an open
neighborhood of $0 \subset W$ invariant under the standard
$U(2)$-action.


\section{Trianalytic subvarieties of compact tori: an
erratum}\label{erratum} 


In our previous paper \cite{kv}, a grave error was found. In the
paragraph 2.3 (page 459 of \cite{kv}), we say ``The right-hand side
obviously depends continuously on the point $a \in F_\Delta$''. This
is, unfortunately, not true. So, even if the rest of the arguments is
valid, the proof of the main theorem is wrong. As we show in the
Appendix, the main theorem of \cite{kv} is false: a generic
deformation of the Hilbert scheme of a compact $2$-dimensional
complex torus contains non-trivial trianalytic subvarieties.

However, while Section 7 of \cite{kv} is false, Sections 1-6 are
correct. In this erratum we use the correct results of \cite{kv} and
some additional arguments to prove a weaker version of the main
theorem of \cite{kv} (Theorem~\ref{trian}). Before we formulate this
result, we need to recall some facts from hyperk\"ahler geometry.

For more details on the false Theorem of \cite{kv} and its
counterexample, see Remark~\ref{counter} below and the Appendix to
this paper.

\subsection{Preliminaries on hyperk\"ahler geometry}
\label{_Prelimi_Subsection_}

\begin{defn}
A {\em hyperk\"ahler manifold} is a Riemannian manifold $M$
equip\-ped with a unitary action $\tau:\h \to \End T(M)$ of the
algebra $\h$ of the quaternions in the tangent bundle $T(M)$ such
that the action $\tau$ is parallel with respect to the Levi-Civita
connection.
\end{defn}

Here {\em unitary} means that for every quaternion $a \in \h$ and
every two tangent vectors $t_1,t_2$ we have 
$$ 
(\tau(a)t_1,t_2) = (t_1,\tau(\overline{a})t_2), 
$$ 
where $\overline{a} \in \h$ is the quaternion, conjugate to $a$.

\hfill

Let $X$ be a hyperk\"ahler manifold. The quaternionic action on the
tangent bundle to $X$ defines an action of the Lie algebra $\su(2)$
on the de Rham complex of the manifold $X$. By \cite{v.su2} this
action commutes with the Laplacian and defines therefore a canonical
$\su(2)$-action on the cohomology spaces $H^\idot(X,\C)$.

Every algebra embedding $I:\C \hookrightarrow \h$ defines by restriction an
almost complex structure $X_I$ on the manifold $X$. Since the almost
complex structure $X_I$ is parallel with respect to the Levi-Civita
metric, it is integrable, and the metric on $X$ is K\"ahler with
respect to the complex structure $X_I$.  The complex structure $X_I$
will be referred to as the {\em induced complex structure} on $X$
associated to the embedding $I$. The corresponding K\"ahler form
will be denoted by $\omega_I$. The set of algebra emebddings 
$\C \hookrightarrow\h$ 
can be naturally identified with the complex projective line $\C
P^1$.

\begin{defn}\label{_triana_Definition_}
A closed subset $Y \subset X$ in a hyperk\"ahler manifold $X$ is called
{\em trianalytic} if it is complex analytic for every induced complex
structure $X_I$ on $X$.
\end{defn}

It is proved in \cite{v.trian} that a subset $Y \subset X_I$
analytic in an induced complex structure $X_I$ is trianalytic if and
only if its fundamental cohomology class is invariant under the
canonical $\su(2)$-action. 

This motivates the following definition. 

\begin{defn}
Let $X$ be a compact complex manifold admitting a hy\-per\-k\"ah\-ler
structure $\HH$. We say that $X$ is 
{\em generic} with respect to $\HH$ if all elements of the group
$$
\bigoplus\limits_p H^{p,p}(X)\cap H^{2p}(X,\Z) \subset H^*(X)
$$
are $SU(2)$-invariant. 
\end{defn}

For every hyperk\"ahler structure, there is at most a countable set
of induced complex structures which are not generic (see, e.g.,
\cite{v.trian}).

\begin{defn}\label{_generic_manifolds_Definition_}
Let $X$ be a compact hyperk\"ahler manifold. An induced
complex structure $I$ on $X$ is called {\em Mumford-Tate generic}
with respect to the hyperk\"ahler structure if for all $n > 0$,
the complex manifold $(X_I)^n$ is generic with respect to the 
hyperk\"ahler structure.
\end{defn}

Fix the standard basis $i,j,k \in \h$ in the space of imaginary
quaternions and consider the associated algebra embeddings $I,J,K:\C
\to \h$. It is easy to check that the form $\Omega = \omega_J +
\sqrt{-1}\omega_K$ is of Hodge type $(2,0)$ with respect to the
complex structure $X_I$ on $X$. Since the from $\Omega$ is closed,
it is holomorphic. Thus every hyperk\"ahler manifold $X$ is
canonically holomorphically symplectic.

For compact manifolds, the converse is also true. Namely, we
have the following corollary of Yau's proof of Calabi conjecture
(\cite{_Yau:Calabi-Yau_}). 

\begin{theorem}[\cite{_Beauville_,_Besse:Einst_Manifo_}] 
\label{_symplectic_=>_hyperkahler_Proposition_}\label{cy}
Let $\langle X, \Omega \rangle$ be a compact complex manifold
equipped with a holomorphic symplectic form $\Omega$, and let
$\omega$ be an arbitrary K\"ahler form on $X$. Then there exists a
unique hyperk\"ahler metric on $X$
with the same K\"ahler class as $\omega$. \endproof
\end{theorem}

Every complex torus $T$ of even dimension is holomorphically symplectic
and admits a K\"ahler metric. Therefore every such torus admits a
hyperk\"ahler structure. 

Let $X$ be a compact hyperk\"ahler manifold of complex dimension $2$
(for example, a $2$-dimensional complex torus). Then the Hilbert
scheme $X^{[k]}$ of $0$-dimensional subschemes of length $k$ in $X$
is a smooth compact complex manifold. It is a manifold of K\"ahler
type, and it carries a canonical holomorphic symplectic form
$\Omega$ (see \cite{Nakajima.Hilbert}). Therefore $X^{[k]}$ admits a
hyperk\"ahler structure (non-canonical, since there is no canonical
K\"ahler metric on $X^{[k]}$).

We can now formulate our main result. Let $T = W/L$ be a
$2$-dimensional complex torus equipped with a holomorphic symplectic
form $\Omega$ and a hyperk\"ahler metric compatible with
$\Omega$. Assume that the lattice $L \subset W$ in the standard
$2$-dimensional complex vector space $W$ is generic in the sense of
Definition~\ref{generic}. Assume also that the torus $T$ is
Mumford-Tate generic with respect to some hyperk\"ahler structure
(Definition~\ref{_generic_manifolds_Definition_}). We say that a
submanifold $Y \subset X$ in a complex manifold $X$ is {\em rigid}
if it admits no non-trivial deformations within $Y$.

\begin{theorem}\label{trian}
Consider a hyperk\"ahler structure on the Hilbert scheme $T^{[k]}$
compatible with the cannical holomorphic symplectic form. Let $X
\subset T^{[k]}$ be an irreducible rigid trianalytic submanifold,
and let $\pi(X) \subset T^{(k)}$ be its image under the canonical
projection $T^{[k]} \to T^{(k)}$.  
\begin{enumerate}
\item The projection $\pi:X \to \pi(X)$ is one-to-one over the
subset of non-singular points in $\pi(X) \subset T^{(k)}$.
\item The complex variety $\pi(X)$ is isomorphic to the quotient
$T'/G$ of a complex torus $T'$ by a finite group $G$. The torus $T'$
is isogenic to a subtorus in the power $T^k$ of the torus $T$. The
group $G$ is generated by real reflexions in a complex vector space
$V_0$, and the $G$-module $V = \Gamma(T',\T(T))$ is isomorphic to
the sum $V = V_0 \oplus V_0$ of two copies of the $G$-module $V_0$. 
\end{enumerate}
\end{theorem}

Note that by virtue of Section 3 of \cite{kv}, for every trianalytic
submanifold $X \subset T^{[k]}$ there exists a hyperk\"ahler
manifold $S$ and a rigid trianalytic submanifold $\wt{X} \subset
T^{[k]}$ whose normalization is isometric to $X \times S$. Therefore
the study of trianalytic submanifolds of $T^{[k]}$ essentially
reduces to the study of rigid trianalytic submanifolds.

Complex tori satisfying conditions of Theorem~\ref{trian} are indeed
generic in the usual sense, that is, they form a dense subset in the
appropriate moduli space, and the complement to this subset is a
countable union of closed analytic subsets of codimension 1
(\cite{v.trian}).

\begin{remark}\label{counter}
The (false) main theorem of \cite{kv} claimed that the only proper
trianalytic subvarieties of the Hilbert scheme $T^{[k]}$ are the
fibers of the summation map $T^{[k]} \to T$. A very simple
counterexample to this statement is constructed as follows. Let
$\iota:T \to T$ be the involution $t \mapsto -t$ of the torus
$T$. Consider the induced involution $\iota:T^{[k]} \to T^{[k]}$ of
the Hilbert scheme $T^{[k]}$. 

The fixed point set $\left(T^{[k]}\right)^\iota \subset T^{[k]}$ is
a union of smooth connected components. Moreover, the involution
$\iota$ preserves the holomorphic symplectic form on $T^{[k]}$ and
the cohomology class of the K\"ahler form. Hence by the uniqueness
statement of Theorem~\ref{cy} it preserves the hyperk\"ahler
structure on $T^{[k]}$.  Therefore each of the connected components
of the fixed point set $\left(T^{[k]}\right)^\iota \subset T^{[k]}$
is a trianalytic subvariety in $T^{[k]}$.

See the Appendix to this paper for precise formulations, proofs and
extensions.
\end{remark}

\subsection{Stratification by diagonals}

Before we begin the proof of Theorem~\ref{trian}, we need to recall
some facts on the stratification by diagonals of the Hilbert scheme
$M^{[k]}$ of a smooth complex manifold $M$ (this was the subject of
Sections~4 and 5 of \cite{kv}).

Let $\Delta = \{a_1,\ldots,a_l\}, \sum a_i = k$ be a partition of
the integer $k$. Denote by $M^{[k]}_\Delta \subset M^{[k]}$ the
subset of points $m \in M^{[k]}$ such that the associated
$0$-dimensional subscheme $\ZZ_m \subset M$ in $M$ is supported on a
finite subset $\{m_1,\ldots,m_l\} \subset M$ of $l$ distinct points
in $M$, and the part of $\ZZ_m$ which is supported on $t_i$ has
length $a_i$. Moreover, denote by $M^{(k)}_\Delta \subset M^{(k)}$
the subset of points $m \in M^{(k)}$ such that the associated subset
of $k$ unordered points in $M$ consists of $l$ distinct points in
$M$, and the point $m_i$ appears with multiplicity $a_i$.

The subsets $M^{[k]}_\Delta \subset M^{[k]}$ and $M^{(k)}_\Delta
\subset M^{(k)}$ are locally closed and form stratifications of
varieties $M^{[k]}$ and $M^{(k)}$. 
Moreover, we have 
$$
M^{[k]}_\Delta = \pi^{-1}\left(M^{(k)}_\Delta\right) \subset M^{[k]},
$$
where $\pi:M^{[k]} \to M^{(k)}$ is the canonical projection. 

A stratum $M^{[k]}_\Delta$ lies in the closure of stratum
$M^{[k]}_{\Delta'}$ if and only if the partition
$\Delta'=\{b_1,\ldots,b_n\}$ of the integer $k$ is a subdivision of
the partition $\Delta = \{a_1,\ldots,a_l\}$. In other words, we must
have $\Delta' = \Delta_1 \cup \cdots \cup \Delta_l$, where
$\Delta_i$ is a partition of the integer $a_i$. Assume that this is
the case. Then for every point $m \subset M^{[k]}_\Delta$ and a
local chart $U^{[\Delta]} \subset M^{[k]}$ near the point $m$, the
intersection $U^{[\Delta]}_{\Delta'} = U^{[\Delta]} \cap
M^{[k]}_{\Delta'}$ decomposes into the direct product
\begin{equation}\label{subdivision}
U^{[\Delta]}_{\Delta'} = U^{[a_i]}_{\Delta_1} \times \cdots \times
U^{[a_l]}_{\Delta_l}.
\end{equation}
An analogous decomposition holds for the local chart $U^{(\Delta)}$
in the symmetric power $M^{(k)}$ near the point $\pi(m) \subset
M^{(k)}$. 

For a partition $\Delta = \{a_1,\ldots,a_l\}$ of the integer $k$,
let $\Sigma_\Delta$ be the group of transpositions of $l$ letters
which preserve the numbers $a_1,\ldots,a_l$ (if all the numbers
$a_1,\ldots,a_l$ are distinct, then the group $\Sigma_\Delta$ is
trivial). Consider the product $M^l = M_1 \times \cdots \times M_l$
of $l$ copies of the manifold $M$. The group $\Sigma_\Delta$ acts on
$M^l$ by transpositions, and the stratum $M^{(k)}_\Delta \subset
M^{(k)}$ is isomorphic to the quotient
$$
M^{(k)}_\Delta \cong \left(M_1 \times \cdots \times M_l \setminus
\Diag\right)/\Sigma_\Delta,
$$
where $\Diag \subset M^l$ is the subset of diagonals. We will denote
by
$$
\wt{M}^{[k]}_\Delta = \left(M^l \setminus \Diag \right)
\times_{M^{(k)}_\Delta} M^{[k]}_\Delta
$$ 
the associated $\Sigma_\Delta$-cover of the stratum $M^{[k]}_\Delta
\subset M^{[k]}$. Points in the variety $\wt{M}^{[k]}_\Delta$
correspond to pairs of
\begin{enumerate}
\item a point $\langle m_1,\ldots,m_l \rangle \in M^l \setminus \Diag$,
and 
\item a $0$-dimensional subscheme $\ZZ \subset M$ supported on the
subset
$$
\{m_1,\ldots,m_l\} \subset M
$$ 
such that the part of $\ZZ$ supported on $m_i$ has length $a_i$.
\end{enumerate}
The induced projection $\pi:\wt{M}^{[k]}_\Delta \to M^l \setminus
\Diag$ sends such a pair to the point
$$
\langle m_1,\ldots,m_l\rangle \in M^l \setminus \Diag, 
$$
and the fiber of the projection $\pi$ over a point $\langle
m_1,\ldots,m_l\rangle \in M^l$ is canonically isomomorphic to the
product
$$
F_{m_1,\ldots,m_l} = F_{m_1} \times \cdot \times F_{m_l},
$$
where $F_{m_i}$ is the punctual Hilbert scheme of $0$-dimensional
subschemes in $M$ of length $a_i$ supported at the point $m_i$. 

Assume now that the holomorphic tangent bundle $\T(M)$ is trivial,
and that we are given a commuting frame in $\T(M)$, that is, a
subspace of commuting vector fields $V \subset \Gamma(M,\T(M))$
which freely generate the bundle $\T(M)$ (this happens naturally in
all our examples of a torus $M=T$, a vector space $M=V$ and an open
neighborhood of $0$ in a vector space $M=U \subset V$). Then the
punctual Hilbert scheme $F_{m_i}$ can be identified canonically with
the punctual Hilbert scheme $F_{a_i}$ of $0$-dimensional subschemes
of length $a_i$ in $V$ supported at $0 \subset V$, and we have the
following.

\begin{lemma}[\cite{kv}, Lemma~5.5]\label{triv}
There exists a canonical splitting
\begin{equation}\label{product}
\wt{M}^{[k]}_\Delta = \left( M^l \setminus \Diag \right) \times
F_\Delta
\end{equation}
of the variety $\wt{M}^{[k]}_\Delta$ into a direct product of the
open subset $M^l \setminus \Diag \subset M^l$ and the product
$$
F_\Delta = F_{a_1} \times \cdots \times F_{a_l}
$$ 
punctual Hilbert schemes of subschemes in $V$ of lengths
$a_1,\ldots,a_l$ supported at $0 \subset V$. Consequently, the
stratum $M^{[k]}_\Delta$ is canonically isomorphic to the quotient
$$
M^{[k]}_\Delta = \left(\left( M^l \setminus \Diag \right) \times
F_\Delta\right) / \Sigma_\Delta,
$$
where $\Sigma_\Delta$ acts on $M^l$ and on $F_\Delta$ by
transpositions. 
\end{lemma}

(In \cite{kv} this is formulated only for a torus $T$, and the proof
uses the abelian group structure on the torus $T$. It is easy to see
that it is enough to have a commuting holomorphic frame $V \subset
\Gamma(M,\T(M))$.) 

The splitting \eqref{product} is functorial with respect to open
embeddings which preserve the commuting holomorphic frame. In
particular, consider a point $m \in M^{[k]}_\Delta$ in a stratum
$M^{[k]}_\Delta \subset M^{[k]}$ which lies in the closure of a
larger stratum $M^{[k]}_{\Delta'} \subset M^{[k]}$. Assume given a
{\em flat} local chart $U^{[\Delta]} \subset U^{[k]}$ near the point
$m$. Consider the intersection $U^{[\Delta]}_{\Delta'} =
U^{[\Delta]} \cap M^{[k]}_{\Delta'}$, and let
$$
U^{[\Delta]}_{\Delta'} = U^{[a_i]}_{\Delta_1} \times \cdots \times
U^{[a_l]}_{\Delta_l}
$$
be the direct product decomposition \eqref{subdivision}. Then the
splitting \eqref{product} for the stratum $M^{[k]}_{\Delta'}$
induces a splitting for the intersection $U^{[\Delta]}_{\Delta'}$,
and this induced splitting coincides with the direct product of
splittings associated to strata $U^{[a_i]}_{\Delta_i}$ in the
Hilbert schemes $U^{[a_i]}$.

Furthermore, assume that $M = T = W/L$ is a $2$-dimensional torus,
and that $U \subset W$ is a $U(2)$-invariant neighborhood of $0
\subset W$. Then the splitting \eqref{product} is compatible with
the canonical $U(2)$-action on the Hilbert schemes $U^{[a_i]}$ and
on the local chart $U^{[\Delta]}$. In other words, the
$U(2)$-action also splits into a product of two commuting
$U(2)$-actions. The first one, which we will call {\em horizontal},
is induced by the standard $U(2)$-action on the local chart
$U^{(\Delta)} \subset T^{(k)}$. The second one, which we will call
{\em vertical}, is induced by the action on the punctual Hilbert
scheme $F_\Delta$.

In terms of the modular data (i)-(ii) above, the horizontal action
affects only the points $m_1,\ldots,m_l \in U$, while the vertical
action does not affect $m_1,\ldots,m_l$, and only changes the
$0$-dimensional subschemes supported at these points.

Note that the vertical $U(2)$-action is induced by an action which
is defined on the whole stratum $T^{[k]}_{\Delta'}$ and which does
not depend the choice of a smaller stratum $T^{[k]}_\Delta$, a point
$m \in M^{[k]}_\Delta = T^{[k]}_\Delta$ and the local chart
$U^{[\Delta]}$ in $T^{[k]}$ near the point $m$.

\subsection{Proof of Theorem~\ref{trian}}

We can now begin the proof Theorem~\ref{trian}. First we formulate
the collection of results of \cite{kv} which we will need for the
proof. Let $T$ be a $2$-dimensional complex torus which is Mumford-Tate generic
with respect to some hyperk\"ahler structure.  Say that a closed
subvariety $X \subset T^{[k]}$ in the Hilbert scheme $T^{[k]}$ {\em
lies generically} in a stratum $T^{[k]}_\Delta \subset T^{[k]}$ iff
$X$ lies in the closure of the stratum $T^{[k]}_\Delta$, while the
intersection $X \cap T^{[k]}_\Delta$ is dense in
$X$. Every irreducible subvariety in $T^{[k]}$ lies generically in
one and only one stratum $T^{[k]}_\Delta$.

\begin{prop}[\cite{kv}, Section~6]\label{old}
Let $X \subset T^{[k]}$ be an irreducible trianalytic subvariety in
the Hilbert scheme of the torus $T$. Assume that $X$ lies
generically in the stratum $T^{[k]}_\Delta \subset T^{[k]}$. Let
$\pi(X) \subset T^{(k)}$ be the image of the subvariety $X$ under
the projection $\pi:T^{[k]} \to T^{(k)}$. Moreover, let $\wt{X}
\subset \wt{T}^{[k]}_\Delta$ be the preimage of the subset $X \cap
T^{[k]}_\Delta \subset T^{[k]}_\Delta$ under the canonical Galois
covering $\wt{T}^{[k]}_\Delta \to T^{[k]}_\Delta$.
\begin{enumerate}
\item the subvariety $\pi(X) \subset T^{(k)}$ is of the form $T'/G$,
where $T' = V'/L' \subset T^k$ is a subtorus, $V' \subset W^k$ is
the associated linear subspace, $L' = V' \cap L^k$ is the induced
lattice in $V'$, and $G = \Norm T' \subset S_k$ is the normalizer
subgroup of the subtorus $T'$. 
\item The subspace $V' \subset W^k$ is
invariant under the canonical $U(2)$-action on the vector space
$W^k$. 
\item The projection $X \to \pi(X)$ is finite and \'etale over an
open dense subset in $\pi(X)$. Moreover, the image of the subset
$\wt{X} \subset \wt{T}^{[k]}_\Delta$ under the natural projection
$$
\wt{T}^{[k]}_\Delta = \left(T^l \setminus \Diag\right) \times
F_\Delta \to F_\Delta
$$
is a finite subset in $F_\Delta$. \endproof
\end{enumerate}
\end{prop}

The new result which we need for the proof of Theorem~\ref{trian} is
the following. Let $X \subset T^{[k]}$ be an irreducible rigid trianalytic
subvariety in the Hilbert scheme $T^{[k]}$ of a Mumford-Tate generic
$2$-dimensional complex torus $T$. Let $U \subset W$ be a
sufficiently small $U(2)$-invariant open neighborhood of $0 \in W$,
and consider the associated flat local chart $U^{[\Delta]} \subset
T^{[k]}$ near a point $x \in X \cap T^{[k]}_\Delta$ for some stratum
$T^{[k]}_\Delta \subset T^{[k]}$. 

\begin{prop}\label{u2.inv}
The intersection $X \cap U^{[\Delta]} \subset U^{[\Delta]}$ is
invariant under the canonical $U(2)$-action on 
$U^{[\Delta]}\subset W^{[\Delta]}$.
\end{prop}

\proof Let $\U$ be an arbitrary vector in the Lie algebra $\U(2)$ of
the Lie group $U(2)$. For every point $x \in X \cap T^{[k]}_\Delta$,
the $U(2)$-action on the local chart $U^{[\Delta]}$ near $x$ induces
an action of the Lie algebra $\U(2)$. Let $\TT_x$ be the holomorphic
vector field $U^{[\Delta]}$ associated to $\U \subset \U(2)$. After
restricting to the intersection $X \cap U^{[\Delta]}$, it defines,
in turn, a local holomorphic section $\n_x$ of the normal bundle
$\N(X)$ to $X$ in $T^{[k]}$. We begin by proving that all these
local sections glue together and define a global section of the
bundle $\N(X)$.

Assume that the irreducible subvariety $X$ lies generically in a
stratum $T^{[k]}_{\Delta'}$ for a partition $\Delta' =
\{a_1',\ldots,a_n'\}$. Since the intersection $X \cap
T^{[k]}_{\Delta'} \subset X$ is dense, it suffices to prove that the
restrictions of the local sections $\n_x$ to the intersections $X_0
= X \cap T^{[k]}_{\Delta'} \cap U^{[\Delta]}$ taken for different
points $x$ glue together to give a global holomorphic section in
$\Gamma(X \cap T^{[k]}_{\Delta'},\N(X))$.

Fix a point $x \in X \cap T^{[k]}_\Delta$ and let $\TT_x = \TT_{hor}
+ \TT_{vert}$ be the decomposition of the vector field $\TT_x$ into
the parts corresponding to the horizontal and the vertical action of
the group $U(2)$ on $U^{[\Delta]}_{\Delta'} = T^{[k]}_{\Delta'} \cap
U^{[\Delta]}$. By Proposition~\ref{old}, the variety $X_0$ lies in a
finite set of fibers of the projection onto the second factor in the
splitting \eqref{product}. Moreover, the subvariety $\pi(X_0) = V
\cap U^n \subset W^n$ is invariant under the standard $U(2)$-action
on $W^n$. Therefore the subvariety $X_0 \subset
U^{[\Delta]}_{\Delta'}$ is preserved by the horizontal
$U(2)$-action.

Consequently, the vector field $\TT_{hor}$ is tangent to the
subvariety $X_0$. (More precisely, the section $\TT_{hor}$ of the
tangent bundle $\T(T^{[k]})$ becomes a section of the holomorphic
tangent bundle $\T(X_0) \subset \T(T^{[k]})$ after the restriction
to $X_0 \subset T^{[k]}$.) Therefore the vector fields $\TT_x$ and
$\TT_{vert}$ define the same section of the normal bundle $\N(X)$.

However, the vertical action of $U(2)$ on $T^{[k]}_{\Delta'} \cap
U^{[\Delta]}$ can be defined globally on the whole
$T^{[k]}_{\Delta'}$ and does not depend on the choice of the point
$x \in x$ and the local chart $U^{[\Delta]}$. Therefore the vector
fields $\TT_{vert}$ for different points $x$ glue together to
produce a global section of the tangent bundle $\T(T^{[k]}$ to
$T^{[k]}$ over $X \cap T^{[k]}_{\Delta'}$. This proves that the
sections $\n_x$ of the normal bundle $\N(X)$ also glue together to a
global section over $X \cap T^{[k]}_{\Delta'}$, hence also over the
whole $X$.

But by \cite{v.trian} (see also Section~3 of \cite{kv}), all
trianalytic subvarieties are unobstructed, that is, every global
section of a normal bundle $\N(X)$ to a trinalytic submanifold $X
\subset T^{[k]}$ defined a local deformation of $X$. Since by
assumption $X$ is rigid, all global sections of $\N(X)$ must
vanish. Therefore all the local sections $\n_x$ also vanish. Since
the vector $\U \in \U(2)$ is arbitrary, this implies that for every
local chart $U^{[\Delta]}$ the intersection $X \cap U^{[\Delta]}$ is
$U(2)$-invariant.  \endproof

We can now prove Theorem~\ref{trian}. 

\proof[Proof of Theorem~\ref{trian}.] 
Assume that the irreducible rigid subvariety $X \subset T^{[k]}$
lies generically in a stratum $T^{[k]}_\Delta$, in other words, the
intersection $X_0 = X \cap T^{[k]}_\Delta \subset X$ is dense in
$X$. Let $x \in X_0$ be an arbitrary point. By
Proposition~\ref{u2.inv} the subvariety $X_0$ is invariant under the
canonical $U(2)$-action corresponding to the local chart at $x$. In
particular, the point $x \in F_\Delta = \pi^{-1}(\pi(x))$ itself is
$U(2)$-invariant. This means that the corresponding $0$-dimensional
subscheme in $W$ supported at $0 \in W$ is $U(2)$-invariant. But it
is easy to see that for every $n$ there exists at most one such
subscheme of length $n$ (\cite{_Verbitsky:Hilbert_}, Sublemma 5.7). 
Therefore the projection $\pi:X_0 \to
\pi(X_0)$ is one-to-one.

To prove \thetag{i}, it remains to prove that $\pi:X \to \pi(X)$ is
one-to-one over the nonsingular part $U \subset \pi(X)$ of $\pi(X) =
T'/G$. Denote by $\wh{U} = \pi^{-1}(U) \cap X$. Since we know that
it is one-to-one over a dense subset $X_0 \subset X$, it suffices to
prove that $\pi:\wh{U} = \pi^{-1}(U) \cap X \to U$ is \'etale.

To do this, consider the symplectic form on the $2$-dimensional
vector space $W$ associated to the hyperk\"ahler structure on the
torus $T$. This form induces an $S_k$-invariant symplectic form
$\Omega_0$ on the vector space $W^k$ and the holomorphic symplectic
form $\Omega$ on the Hilbert scheme $T^{[k]}$. By Lemma~6.2 of
\cite{kv} the restriction of the form $\Omega_0$ to the vector
subspace $V \subset W^k$ associated to the subtorus $T'$ is
non-degenerate. Therefore it induces a holomorphic symplectic form
$\Omega_0$ on the non-singular part $U \subset T'/G$. By
\cite[Proposition 4.5]{_Verbitsky:Hilbert_}, 
the pullback $\pi^*(\Omega_0)$
coincides with the restriction of the form $\Omega$ on $T^{[k]}$ to
the subset $\wh{U} \subset T^{[k]}$.

Since $X \subset T^{[k]}$ is trianalytic, this implies that the
pullback $\pi^*(\Omega_0)$ is non-degenerate. But both $U$ and
$\wh{U}$ are smooth. Therefore the differential $d\pi:\T(\wh{U}) \to
\pi^*\T(U)$ of the map $\pi:\wh{U} \to U$ is an isomorphism, and the
map $\pi:\wh{U} \to U$ is indeed \'etale. This proves \thetag{i}.

To prove \thetag{ii}, we note that by Proposition~\ref{old} the
image $\pi(X) \subset T^{(k)}$ is of the form $T'/G$, where $T'
\subset T^k$ is a subtorus, and $G = \Norm T' \subset S_k$ is the
normalizer subgroup of the subtorus $T' \subset T^k$ in the
symmetric group $S_k$. Together with \thetag{i} this means that $X
\subset T^{[k]}$ is a partial resolution of Hilbert type in the
sense of Definition~\ref{hilbert.type}. In particular, $T' = V/L'$,
where $V = \Gamma(T',\T(T) \subset W^k$ is a linear subspace, and
$L' = L^k \cap V \subset V$ is the induced lattice.

The finite group $G$ carries a canonical three-step filtration
$$
G_0 \subset G_1 \subset G,
$$
where $G_0 \subset G$ is the subgroup of elements which act
trivially on the torus $T'$, and $G_1 \subset G \subset S^k$ is the
subgroup of elements which act trivially on the vector subspace $V
\subset W^k$. The action of the group $G$ on the torus $T'$ factors
through the quotient $\wh{G} = G/G_0$, and we have $\wh{G}_1 =
G_1/G_0 \subset \wh{G}$. Denote by $\sigma:\wh{G} \to
\wh{G}/\wh{G}_0$ the quotient map.

The subgroup $\wh{G}_1 = G_1/G_0 \subset \wh{G}$ acts on the torus
$T'$ by translations of finite order. Therefore it is
abelian. Moreover, every element $g \in \wh{G}_1$ acts on the torus
$T'$ without fixed points. 

Conversely, every element $g \in \wh{G}$ which acts on $T'$ without
fixed points must act by a translation. Therefore $g$ acts trivially
on the vector space $V = \Gamma(T',\T(T))$, in other words, belongs
to $\wh{G}_1 \subset \wh{G}$. 

For every point $t \in T'$, let $G_t \subset \wh{G}$ be the
stabilizer subgroup of the point $t$, and let $\sigma(G_t) \subset
\wh{G}/\wh{G}_1$ be its image under the quotient map. Since
$\wh{G}_1 \subset \wh{G}$ is precisely the subgroup of elements
which act on $T'$ without fixed points, the canonical map $G_t \to
\sigma(G_t)$ is an isomorphism, and the quotient
$$
\wh{G}/\wh{G}_1 = \bigcup_{t \in T'} \sigma(G_t)
$$
is the union of the subgroups $\sigma(G_t) = G_t \subset
\wh{G}/\wh{G}_1$ for different points $t \in T'$.

But the quotient $\wh{G}/\wh{G}_1$ acts naturally on the subspace
$V_0 = V^{U(1)} \subset V$ of $U(1)$-invariants in $V$. Moreover,
Theorems~\ref{complex.refl} and \ref{real.refl} imply that the
stabilizer $G_t \in \wh{G}/\wh{G}_0$ of an arbitrary point $t \in
T'$ acts on $V_0$ as a group generated by real reflections. Therefore
the whole group $\wh{G}/\wh{G}_1$ acts on the space $V_0$ as a group
generated by real reflections. 

To finish the proof, it remains to notice that 
$$
\pi(X) = T'/G =
\left(T'/\wh{G}_1\right)/\left(\wh{G}/\wh{G}_1\right),
$$
and the quotient $T'/\wh{G}_1$ is a complex torus isogenic to $T'
\subset T^k$. Replacing $G$ with $\wh{G}/\wh{G}_1$ and $T'$ with
$T'/\wh{G}_1$ proves \thetag{ii}.  \endproof


\section{Appendix. Existence of subvarieties of generalized Kummer
varieties}
\label{_Appe_Exi_Section_}

\subsection{The counterexamples}

Let $T$ be a complex torus of dimension $2$.  Consider the Hilbert
scheme $T^{[n+1]}$ of $n+1$ points on $T$. This is a complex variety
of dimension $2(n+1)$. The commutative group structure on the torus
$T$ defines a summation map $\Sigma:T^{n+1} \to T$, which induces a
morphism $\Sigma:T^{[n+1]} \to T$. It is easy to see
that $\Sigma$ coincides with the Albanese map.

\begin{defn}\label{_generali_Kummer_Definition_}
The {\em generalized Kummer variety} $K^{[n]}$ associated to the
torus $T$ is the preimage $\Sigma^{-1}(0) \subset T^{[n+1]}$ of the
zero $0 \in T$.
\end{defn}

In \cite{kv}, the following false theorems were stated.

\hfill

{\bf Theorem (false).} {\it
Let $T$ be a $2$-dimensional complex torus equipped with a
hyperk\"ahler structure. Assume that the complex structure on $T$ is
Mumford-Tate generic, and consider a generalized Kummer variety
$K^{[n]}$ associated to $T$. 

For any hyperk\"ahler structure on $K^{[n]}$ compatible with the
canonical holomorphic $2$-form, every irreducible trianalytic
subvariety $X \subset K^{[n]}$ is either the whole $K^{[n]}$ or a
single point.}

\hfill

{\bf Theorem (false).} {\it Assume that the Kummer variety $K^{[n]}$
is equipped with a complex structure $I$, and assume that $I$ is
Mumford-Tate generic with respect to some hy\-per\-k\"ah\-ler
structure on $K^{[n]}$ compatible with the canonical holomorphic
$2$-form.

Then every irreducible analytic subvariety $X \subset K^{[n]}_I$ is
either the whole $K^{[n]}$ or a single point. 
}

\hfill

Here we offer the counterexamples to these statements.

Notice that by a ``hyperk\"ahler structure on a complex manifold''
we understand a hyperk\"ahler structure inducing the complex structure.
A manifold of K\"ahler type is a complex manifold admitting a K\"ahler 
metric.

\begin{theorem}\label{_Involution_examples_Theorem_}
Let $T$  be a compact torus, $K^{[n]}$ the $n$-th generalized Kummer
variety of $T$, $n>1$, and $X$ be a manifold of K\"ahler type
which is a deformation of $K^{[n]}$. Then $X$ admits a non-trivial 
complex analytic involution $\iota:\; X\arrow X$,
and the fixed set of $\iota$ has positive dimension.
Moreover, for any hyperk\"ahler structure $\HH$ on $X$,
the involution $\iota$ is compatible with $\HH$.
\end{theorem}

The proof of Theorem~\ref{_Involution_examples_Theorem_}
is given in Subsection~\ref{_example_Subsection_}.

\begin{remark}
The fixed set of an involution is always a union of smooth
manifolds. From Theorem~\ref{_Involution_examples_Theorem_},
it follows that the fixed point set of the involution $\iota$
is always trianalytic. Thus,
Theorem~\ref{_Involution_examples_Theorem_}
gives counterexamples to both false theorems of 
\cite{kv}.
\end{remark}

The construction of the involution $\iota$ was suggested by
A. Kuznetsov for the Hilbert scheme of $\C^2$. 
Let $\iota_0:\; \C^2\arrow \C^2$ be the involution
mapping $(x,y) \arrow (-x,-y)$. By functoriality,
$\iota_0$ is naturally extended to the Hilbert scheme
$(\C^2)^{[n]}$ of $\C ^2$. A. Kuznetsov conjectured that
this involution is compatible with the natural hyperk\"ahler
structure on $(\C^2)^{[n]}$, constructed by 
Kronheimer and Nakajima (see, e. g. \cite{Nakajima.Hilbert}).
This conjecture was not proven, though the proof seems to be
straightforward. However, using Calabi-Yau
(Theorem~\ref{_symplectic_=>_hyperkahler_Proposition_}),
it is possible to prove that the analogous involution
is compatible with the hyperk\"ahler structure
on the Hilbert scheme of a torus.

\subsection[Twistor paths and diffeomorphisms preserving
trianalytic \protect\\subvarieties]{Twistor paths and
diffeomorphisms preserving trianalytic subvarieties}

Let $M$ be a hyperk\"ahler manifold. Consider the product manifold 
$X = M\times S^2$. Embed the sphere $S^2 \subset {\Bbb H}$ 
into the quaternion algebra ${\Bbb H}$ 
as the subset of all quaternions $J$ with $J^2 = -1$. For every point
$x = m \times J \in X = M \times S^2$ the tangent space $T_xX$ is
(using Levi-Civita connection) 
decomposed as $T_xX = T_mM \oplus T_JS^2$. Identify $S^2 = \C P^1$
and let $I_J:T_JS^2 \to T_JS^2$ be the complex structure operator. Let
$I_m:T_mM \to T_mM$ be the complex structure on $M$ induced by $J \in S^2
\subset {\Bbb H}$.

The operator $I_x = I_m \oplus I_J:T_xX \to T_xX$ satisfies $I_x \circ I_x =
-1$. It depends smoothly on the point $x$, hence defines an almost complex
structure on $X$. This almost complex structure is known to be integrable
(see \cite{_Salamon_}). The complex manifold $X$ is called 
{\bf the twistor space of} $M$.

The twistor space is fibered over $\C P^1$, and the corresponding
map $\pi:\; X\arrow \C P^1$ is called {\bf the twistior
fibration}. The map $\pi$, considered as a deformation of $M$ over
the base $\C P^1$, gives an embedding of $\C P^1$ into the moduli
$Comp$ of deformations of complex structures on $M$.  Such an
embedding is called {\bf the twistor curve, associated with the
hyperk\"ahler structure $\HH$}.  A sequence of connected twistor
curves is called {\bf a twistor path}. In \cite{_coho_announce_} it
was proven that any two points in the moduli space $Comp$ can be
connected by a twistor path.  For more detailed definitions and
results concerning the twistor paths, the reader is referred to
\cite{_Verbitsky:examples_} (Subsection 10.1).

\begin{defn}
Let  $P_0$, ..., $P_n\subset Comp$ be
a sequence of 
twistor curves, supplied with an intersection point
$x_{i+1}\in P_i\cap P_{i+1}$ for each $i$. We say that
$\gamma= P_0, ..., P_n, x_1, ..., x_n$ is 
a {\bf twistor path}. Let $I$, $I'\in Comp$.
We say that $\gamma$ is {\bf a twistor path
connecting $I$ to $I'$} if $I\in P_0$ and $I'\in P_n$.
The lines $P_i$ are called {\bf the edges},
and the points $x_i$ {\bf the vertices}
of a twistor path.
\end{defn}

Recall that in \ref{_generic_manifolds_Definition_},
we defined induced complex structures which are 
generic with respect to a hyperk\"ahler structure.

\begin{defn}
Let $I$, $J\in Comp$ and $\gamma= P_0, ..., P_n$ be a 
twistor path from $I$ to $J$, which corresponds to the 
hyperk\"ahler structures 
$\HH_0$, ..., $\HH_n$. We say that $\gamma$ is {\bf admissible}
(resp. {\bf Mumford-Tate admissible})
if all vertices of $\gamma$ are generic
(resp. Mumford-Tate generic) with respect 
to the corresponding edges.
\end{defn}

Clearly, every twistor path $\gamma$
induces a diffeomorphism \[ \mu_\gamma:\; (M,I)\arrow (M,I'). \]
In \cite{_coho_announce_}, Subsection 5.2,
we studied algebro-geometrical properties of this
diffeomorphism. 

\begin{theorem} \label{_to_triana_Theorem_}
Let $I$, $J\in Comp$ and $\gamma= P_0, ..., P_n$ be an admissible 
twistor path from $I$ to $J$. 
Then $\mu_\gamma$ maps trianalytic subvarieties
to trianalytic subvarieties. 
\end{theorem}

{\bf Proof:} For an induced complex structure of generic type,
every closed complex analytic subvariety is trianalytic.\ Therefore,
diffeomorphism between the generic fibers of 
the twistor fibration,
provided by a single twistor curve, maps complex
subvarieties to complex subvarieties. Taking compositions
of such diffeomorphisms, we obtain a diffeomorphism
mapping complex subvarieties to complex subvarieties, and hence,
trianalytic subvarieties to trianalytic subvarieties.
For more details, see
\cite{_coho_announce_}, Subsection 5.2.
\endproof

\begin{corollary}\label{_Automo_Corollary_}
Let $I$, $J\in Comp$ and $\gamma= P_0, ..., P_n$ be a Mumford-Tate 
admissible twistor path from $I$ to $J$. 
Then $\mu_\gamma$ maps hyperk\"ahler automorphisms
of $(M, I)$ to hyperk\"ahler automorphisms
of $(M, J)$.
\end{corollary}

{\bf Proof:} Taking a graph of a hyperk\"ahler automorphism,
we obtain a trianalytic subvariety in $M\times M$. Applying
Theorem~\ref{_to_triana_Theorem_} to thus subvariety,
we obtain that the diffeomorphism corresponding
to a Mumford-Tate admissible twistor path 
maps hyperk\"ahler automorphisms to 
hyperk\"ahler automorphisms.
\endproof

Notice that a generic deformation of a 
compact hyperk\"ahler manifold 
satisfies $Pic(M_I)=0$.

\begin{theorem}\label{_admi_exi_Theorem_} 
Let $\HH$, $\HH'$ be hyperk\"ahler structures on $M$, and 
$I$, $I'$ be complex structures of general type with respect to 
to and induced by $\HH$, $\HH'$. Assume that $Pic(M_I)=0$.
Then $I$, $I'$ can be connected by an admissible path.
\end{theorem}

{\bf Proof:}  This is \cite{_coho_announce_}, Theorem 5.2.
\endproof

\begin{corollary}\label{_automo_from_gene_Corollary_}
Let $I$ be a complex structure on $M$, such that
$Pic(M,I)=0$, and $\iota$ be an automorphism
of $(M,I)$. Assume that $I$ is induced by 
a hyperk\"ahler structure $\HH$, and $(M, I)$ is
Mumford-Tate generic with respect to $\HH$.
Let $I'$ be a complex structure of K\"ahler type
on $M$ which lies in the same deformation class
as $I$. Then there exists a diffeomorphism
$$\mu:\; (M,I) \arrow (M, I')$$ which
maps $\iota$ to a complex analytic automorphism of $(M, I')$
\end{corollary}
%
%

{\bf Proof:} A diffeomorphism
$\mu$ is constructed from a Mumford-Tate admissible twistor
path, as follows from Theorem~\ref{_admi_exi_Theorem_} 
and Corollary~\ref{_Automo_Corollary_}.
\endproof

\subsection{Examples of subvarieties of generalized Kummer varieties}
\label{_example_Subsection_}

{}From Corollary~\ref{_automo_from_gene_Corollary_},
we obtain the following. Let $\iota$ be an involution of
a compact holomorphic symplectic manifold $(M,I)$, $Pic(M,I)=0$, which
is Mumford-Tate generic with respect to some hyperk\"ahler
structure. Then $\iota$ is compatible with any hyperk\"ahler
structure $\HH$ on $M$ such that $I$ is Mumford-Tate 
generic with respect to $\HH$. Moreover, for any
complex structure $I'$ on $M$ in the same deformation
class, $\iota$ is mapped to a complex automorphism of 
$(M, I')$ by a diffeomorphism, associated with an 
admissible twistor path. Therefore, to prove
Theorem~\ref{_Involution_examples_Theorem_},
it suffuces to prove the following.

\begin{theorem}\label{_exist_invo_gene_Lemma_}
Let $M$, $Pic(M)=0$, 
be a deformation of a generalized Kummer variety,
which is Mumford-Tate generic with respect to some
hyperk\"ahler structure.  Then $M$ admits a non-trivial complex
analytic involution $\iota$. Moreover, the fixed point set of
$\iota$ has positive dimension.
\end{theorem}

{\bf Proof:} 
To prove Theorem~\ref{_exist_invo_gene_Lemma_},
we need to construct an involution $\iota$ on a ``sufficiently
generic'' deformation of $K^{[n]}$. Let $\M$ be the moduli space of 
complex deformations of $M$, and $\M_0\subset \M$
a subset consisting of all complex structures $I'$
on $M$ satisfying the following.
\begin{quote}
{\it The manifold $(M, I')$ admits a hyperk\"ahler structure
$\HH$ inducing a complex structure $I''$ such that
$(M, I'')$ is a generalized Kummer variety for some torus $T$.
}\end{quote}
In other words, $\M_0$ is the set of all deformations 
$(M,I)\in \M_0$  of $K^{[n]}$
which can be connected to a generalized Kummer
variety by a single twistor curve.

\hfill

In the proof of Theorem~\ref{_exist_invo_gene_Lemma_}, 
we use the following lemma.

\begin{lemma} \label{_Open_in_moduli_Lemma_}
In the above assumptions, the set $\M_0\subset \M$ 
is open in $\M$.
\end{lemma} 

{\bf Proof:} The proof is based on an
elementary dimension count. We use the Torelli theorem
which states that the period map from $\M$ to
$\Bbb P H^2(M, \C)$ is etale onto its image, which 
is of codimension one in $\Bbb P H^2(M, \C)$.
The moduli of complex tori is a 
4-dimensional complex space, and the moduli of
deformations of a complex structure on 
a generalized Kummer variety is 5-dimensional.
An extra dimension is given by a choice of
an induced complex structure on a Kummer 
variety.
\endproof

{}From Lemma~\ref{_Open_in_moduli_Lemma_},
it follows that to prove Theorem~\ref{_exist_invo_gene_Lemma_}
it suffices to construct an involution 
$\iota:\; K^{[n]}\arrow K^{[n]}$
that commutes with any hyperk\"ahler structure.
This is done as follows. We have an involution of a torus
$\iota:\; T\arrow T$  mapping $x$ to $-x$. 
By functoriality, this involution is extended 
to an involution of the Hilbert scheme $T^{[n+1]}$.
Since $\iota$ commutes with the Albanese map 
$\Sigma:\;T^{[n+1]}\arrow T$, this map
preserves the generalized Kummer
variety $K^{[n]}\subset T^{[n+1]}$.
Consider the natural projection map
$\pi:\; K^{[n]}\arrow T^{(n+1)}$
from the Kummer variety to the 
symmetric power of $T$. 
Let $C\subset K^{[n]}$ be the singular locus of $\pi$,
and $[C]\in H^2(K^{[n]})$ the corresponding cohomology class. Clearly,
$H^2(K^{[n]})$ is generated by $[C]$ 
and the group $\pi^*(H^2(T^{(n+1)}))\cong H^2(T)$.
Since $\iota$ acts identically on $\pi^*(H^2(T^{(n+1)}))\cong H^2(T)$
and preserves $C$, $\iota$ acts as identity on $H^2(K^{[n]})$.
Therefore, $\iota$ maps a K\"ahler class to itself,
for any K\"ahler class $\omega$ on $K^{[n]}$. By
Theorem~\ref{_symplectic_=>_hyperkahler_Proposition_},
$\iota$ preserves the hyperk\"ahler metric. It is easy
to see that there exists (at most) a unique hyperk\"ahler
structure compatible with a given holomorphic 
symplectic structure and a metric. Therefore, $\iota$
preserves the hyperk\"ahler structure.

We obtained an involution $\iota$
on a ``sufficiently generic'' deformation
of $K^{[n]}$. It remains to show that
this involution has positive-dimensional fixed
set. This is done by an elementary geometric observation.
Assume, for instance, that $n$ is odd, $n+1=2k$.
Then $\iota$ fixes all the $2k$-tuples 
$$
(x_1, y_1, x_2, y_2, ... x_k, y_k)\in T^{(n+1)}
$$
satisfying $x_i=-y_i$. When $x_i, y_i$ are pairwise
distinct, they give a point of the Hilbert scheme
fixed by $\iota$. The closure $D$ of the set of such points
is a component of the fixed point set of $\iota$, and it is obviously
positively-dimensional.  
This proves Theorem~\ref{_exist_invo_gene_Lemma_}.
\endproof

\begin{remark}
It is easy to see that $D$ is birationally equivalent to a Hilbert
scheme of a Kummer K3 surface associated with $T$.
\end{remark}

\hfill

{\bf Acknowledgements:} We are grateful to 
V. Batyrev, V. Ginzburg, A. Kuznetsov, and T. Pantev
for interesting discussions.

\hfill

\end{document}